# SEMIPARAMETRIC ESTIMATION OF A TWO-COMPONENT MIXTURE MODEL

By Laurent Bordes, Stéphane Mottelet and
Pierre Vandekerkhove

*Université de Technologie de Compiègne, Université de Technologie de Compiègne and Université de Marne-la-Vallée*

Suppose that univariate data are drawn from a mixture of two distributions that are equal up to a shift parameter. Such a model is known to be nonidentifiable from a nonparametric viewpoint. However, if we assume that the unknown mixed distribution is symmetric, we obtain the identifiability of this model, which is then defined by four unknown parameters: the mixing proportion, two location parameters and the cumulative distribution function of the symmetric mixed distribution. We propose estimators for these four parameters when no training data is available. Our estimators are shown to be strongly consistent under mild regularity assumptions and their convergence rates are studied. Their finite-sample properties are illustrated by a Monte Carlo study and our method is applied to real data.

**1. Introduction.** Cumulative distribution functions (c.d.f.) of $p$-variate multi-component mixture models are generally defined by

$$(1) \qquad G(x) = \sum_{j=1}^{k} \lambda_j F_j(x), \qquad x \in \mathbb{R}^p,$$

where the unknown mixture proportions $\lambda_j$ ($\lambda_j \geq 0$ and $\sum_{j=1}^{k} \lambda_j = 1$) and the unknown c.d.f. $F_j$ are to be estimated. It is commonly assumed that the $F_j$'s belong to a parametric family, which means that the space of unknown parameters is reduced to a Euclidean set, leading to parametric inference. There is an extensive literature on this subject, including the monographs of Everitt and Hand [16], Titterington, Smith and Makov [40], McLachlan and Basford [28] and McLachlan and Peel [29]. The main types of estimators that have been proposed are the following: maximum likelihood (e.g., [7, 24, 25, 35]), minimum chi-square (e.g., [11]), method of moments (e.g., [26]), Bayesian approaches (e.g., [13, 15]) and techniques based on moment generating functions (e.g., [34]). Note that the number of components $k$ in model (1) may also be an unknown parameter to be estimated, leading to various rates of convergence for maximum likelihood estimators, as discussed by Chen [6]. In this case, the selection of a model is an important topic; see, for example, [10, 22, 23].

The choice of a parametric family for the $F_j$'s may be difficult when little is known about subpopulations. However, models of type (1) are generally nonparametrically nonidentifiable without additional assumptions. This is no longer true when training data are available, that is, when some data are of known origin with respect to the components of the mixture distribution. In this case nonparametric techniques can be applied; see, for example, [4, 17, 21, 31, 33, 36, 39, 40]. As Hall and Zhou [18] state, "very little is known about nonparametric inference in mixtures without training data." These authors looked at $p$-variate data drawn from a mixture of two distributions, each having independent components, and proved that, under mild regularity assumptions, their model is identifiable for $p \geq 3$. They proposed root-$n$ consistent estimators of the $2p$ univariate marginal distributions and the mixing proportion. In a working paper Kitamura [20] investigates identifiability of type (1) models with the presence of covariates.

Note that even if model (1) is not nonparametrically identifiable, there exist, for $p = 1$ and $k = 2$, many real data sets in the statistical literature for which such a model is used under parametric assumptions for the $F_j$'s. For example, Azzalini and Bowman [1] provided data on the length of intervals between eruptions and the duration of the eruption for the Old Faithful Geyser in Yellowstone National Park. Another example deals with average amounts of precipitation (rainfall) in inches for United States cities (from the Statistical Abstract of the United States, 1975; see [30]). These two data sets are available in the R statistical package. Moreover, in some studies, the only parameters of interest are mixture proportions, in which case components $F_j$ in model (1) are nuisance parameters (see, e.g., [8]). In this paper we consider the two-component identifiable restriction of model (1) defined by

$$(2) \qquad G(x) = \lambda F(x - \mu_1) + (1 - \lambda)F(x - \mu_2), \qquad x \in \mathbb{R}.$$

Unknown parameters are the c.d.f. $F$ of a symmetric distribution, two real location parameters $\mu_1$ and $\mu_2$ and the mixing proportion $\lambda$. Note that this model has also been studied by Hunter, Wang and Hettmansperger [19] in an independent work. Model (2) above is called *semiparametric* inasmuch as the unknown parameters can be separated into a functional part $F$ and a Euclidean part $(\mu_1, \mu_2, \lambda)$. Note that such a model should be distinguished from the so-called semi- or nonparametric mixture models (e.g., [27]) where $G$ is defined by

$$(3) \qquad G(x) = \int_{\mathbb{R}} F(x; \theta) \, dH(\theta), \qquad x \in \mathbb{R},$$

where $F$ belongs to a parametric family and $H$ is an unknown distribution function on $\mathbb{R}$. However, as noted by Lindsay and Lesperance [27], there is a link between models (1) and (3) if $H$ is discrete with $k$ points of support. Of course, such a link exists between models (2) and (3) by assuming that, in the latter model, $F(\cdot; \theta) = F(\cdot - \theta)$ with $F$ in the c.d.f. family $\mathcal{F}$ of symmetric distributions, and that $H$ puts masses $\lambda$ and $1 - \lambda$ at points $\mu_1$ and $\mu_2$, respectively.

One of the fundamental issues with mixture models of type (1) is to provide identifiability results. When the $F_j$'s belong to certain specific parametric families (e.g., the continuous exponential family), identifiability results are available; see, for example, [2, 5, 38]. More is needed when we aim to estimate the $F_j$'s nonparametrically (see [18] for the two components case). Working with model (2), we need to prove that $G$ is defined by a unique quadruple $(\lambda, \mu_1, \mu_2, F)$.

The paper is organized as follows. In the next section we give an identifiability result for model (2). In Section 3 we provide a methodology for estimating unknown parameters in our two-component mixture model. Consistency results and convergence rates of our estimators are given in the same section. Our main results are proved in Section 5. In Section 4 finite-sample properties of our estimators are illustrated by a Monte Carlo study and our method is applied to precipitation data.

**2. Identifiability.** First, note that if $F$ in model (2) admits a density function $f$ (an even function), the mixture distribution admits a density function $g$ defined by

$$(4) \qquad g(x) = \lambda f(x - \mu_1) + (1 - \lambda) f(x - \mu_2), \qquad x \in \mathbb{R},$$

where $\theta = (\lambda, \mu_1, \mu_2) \in \Theta = [0, 1/2] \times (\mathbb{R}^2 \backslash \Delta)$ and $\Delta = \{(x, x); x \in \mathbb{R}\}$.

The aim of this section is to investigate identifiability, that is, the possibility of having

$$(5) \qquad \begin{aligned} &\lambda F(x - \mu_1) + (1 - \lambda) F(x - \mu_2) \\ &\quad = \lambda' F'(x - \mu_1') + (1 - \lambda') F'(x - \mu_2') \qquad \forall x \in \mathbb{R}, \end{aligned}$$

for two different quadruples $(\theta, F)$ and $(\theta', F')$ in $\Theta \times \mathcal{F}$, where $\theta' = (\lambda', \mu_1', \mu_2')$ and $\mathcal{F}$ is the c.d.f. set of symmetric distributions. Note that it is sufficient to consider $\lambda \in [0, 1/2]$ because the model is invariant by permutation of $(\lambda, \mu_1)$ and $(1 - \lambda, \mu_2)$. Note also that what we mean by identifiability is not entirely an injectivity condition, since if $\lambda = 0$, we only need to obtain $\lambda' = 0$, $\mu_2 = \mu_2'$ and $F = F'$. Clearly, identifiability fails if we allow $\lambda$ to be equal to $1/2$. Indeed, suppose that $f$ is itself an even mixture density function, for example, $f(x) = h(x - \mu)/2 + h(x + \mu)/2$ with $h$ an even density function. If $g(x) = f(x - \mu_2)$ with $\lambda = 0$, then (5) is obviously satisfied with $\lambda' = 1/2$, $\mu_1' = \mu + \mu_2$, $\mu_2' = \mu_2 - \mu$ and $f' = h$. The main identifiability result is summarized in the following theorem.

THEOREM 2.1. *If* $(\lambda, \mu_1, \mu_2, F)$ *and* $(\lambda', \mu_1', \mu_2', F')$ *are two parameters of* $[0, 1/2) \times (\mathbb{R}^2 \backslash \Delta) \times \mathcal{F}$ *satisfying* (5)*, then* $\lambda = \lambda'$*,* $\mu_2 = \mu_2'$ *and* $F = F'$*, and* $\mu_1 = \mu_1'$ *if* $\lambda > 0$*.*

Hunter, Wang and Hettmansperger [19] have established a similar result for the parametric part $(\lambda, \mu_1, \mu_2)$ of the model. Their results are slightly different from ours since they considered identifiability from the injectivity point of view. They also gave a necessary condition for identifying a type (2) model with three components.

A question which naturally arises concerns the possibility of extending our identifiability result when scale parameters are introduced into model (2). In fact, it is easy to show that such a model is generally not identifiable.

**3. Methodology and theoretical results.** Let $X_1, \ldots, X_n$ be $n$ independent and identically distributed random variables with common c.d.f. $G$ given by model (2). We shall denote by $\theta_0$ and $F_0$ the true values of the unknown Euclidean parameter and the unknown mixed c.d.f. The aim of this section is to propose estimators for $\theta_0$ and $F_0$. Asymptotic results are given with respect to $n \to +\infty$.

The first key idea developed in Section 3.1 is based on the possibility of expressing $F$ as a function of $G$ and $\theta$ (resp. $f$ as a function of $g$ and $\theta$) by inverting the relation (2) [resp. by inverting the relation (4)]. The second key idea, developed in Section 3.2, involves using the symmetry property of $F_0$ in order to propose a contrast function for the Euclidean parameter $\theta$ when $G$ is known. Then, in Sections 3.3 and 3.4, replacing $G$ by the corresponding empirical c.d.f., we propose estimators of $\theta_0$ and $F_0$ and give some asymptotic results for these estimates. These results are obtained under two kinds of assumptions on $F_0$:

C1. $F_0$ is strictly increasing and Lipschitz on $\mathbb{R}$.

C2. $F_0$ is strictly increasing, twice continuously differentiable on $\mathbb{R}$ and $F_0'' \in L^1(\mathbb{R})$.

3.1. *Inversion formula.* Assume that in the mixture model defined by (2) the Euclidean parameter $\theta = (\lambda, \mu_1, \mu_2)$, with $\mu_1 \neq \mu_2$ and $\lambda \in [0, 1/2)$, is known. The key idea consists in rewriting (2) as

$$(6) \qquad F(x) = \frac{1}{1-\lambda} G(x + \mu_2) + \frac{-\lambda}{1-\lambda} F(x + \eta) \qquad \forall x \in \mathbb{R},$$

where $\eta = \mu_2 - \mu_1 \neq 0$, and hence, using (6) as a recurrence formula. Let $\ell$ be a positive integer. By using (6) $\ell$ times, we get

$$(7) \qquad \begin{aligned} F(x) &= \frac{1}{1-\lambda} \sum_{k=0}^{\ell-1} \left( \frac{-\lambda}{1-\lambda} \right)^k G(x + \mu_2 + k\eta) \\ &\quad + \left( \frac{-\lambda}{1-\lambda} \right)^\ell F(x + \ell\eta) \qquad \forall x \in \mathbb{R}. \end{aligned}$$

Let us show that

$$(8) \qquad F(x) = \frac{1}{1-\lambda} \sum_{k \geq 0} \left( \frac{-\lambda}{1-\lambda} \right)^k G(x + \mu_2 + k\eta) \qquad \forall x \in \mathbb{R}.$$

If we denote by $H$ the right-hand side in (8), then by (7) we get, for all $\ell \geq 1$,

$$\|F - H\|_\infty \leq \left( \frac{\lambda}{1-\lambda} \right)^\ell + \frac{1}{1-\lambda} \sum_{k \geq \ell} \left( \frac{\lambda}{1-\lambda} \right)^k$$

$$\leq \left( \frac{\lambda}{1-\lambda} \right)^\ell \left( 1 + \frac{1}{1-2\lambda} \right),$$

where $\|\cdot\|_\infty$ denotes the supremum norm. Since the right-hand side of the above inequality can be made arbitrarily small, it follows that $F = H$. Similarly, working with densities [see (4)] and replacing the supremum norm by the $L^1(\mathbb{R})$-norm $\|\cdot\|_1$, we get

$$(9) \quad f(x) = \frac{1}{1-\lambda} \sum_{k \geq 0} \left( \frac{-\lambda}{1-\lambda} \right)^k g(x + \mu_2 + k\eta) \qquad \text{for } \mu\text{-almost all } x \in \mathbb{R},$$

where $\mu$ is Lebesgue measure on $\mathbb{R}$.

At this point it is convenient to introduce the linear bounded operators $A_\theta$ and $A_\theta^{-1}$ defined by

$$(10) \quad A_\theta = \lambda \tau_{\mu_1} + (1-\lambda) \tau_{\mu_2} \quad \text{and} \quad A_\theta^{-1} = \frac{1}{1-\lambda} \sum_{k \geq 0} \left( \frac{-\lambda}{1-\lambda} \right)^k \tau_{-\mu_2 - k\eta},$$

where $\tau_\mu$ ($\mu \in \mathbb{R}$) is the shift operator defined by $\tau_\mu f = f(\cdot - \mu)$. With the above definitions of $A_\theta$ and $A_\theta^{-1}$, formulae (2) and (4) are equivalent to $G =$ $A_\theta F$ and $g = A_\theta f$, respectively, whereas formulae (8) and (9) are equivalent to $F = A_\theta^{-1} G$ and $f = A_\theta^{-1} g$, respectively.

The interest of the operator $A_\theta^{-1}$ is that if $\theta$ is known, the c.d.f. $F$ may be recovered from a nonparametric estimate $\hat{G}$ of $G$ by considering the reversed estimates $\hat{F} = A_\theta^{-1} \hat{G}$. This also holds for the density $f$, defining $\hat{f} = A_\theta^{-1} \hat{g}$ with $\hat{g}$ a nonparametric estimator of $g$. Unfortunately, the Euclidean parameter $\theta$ is generally unknown and thus we need to propose an estimate of $\theta$ separately. It should be noted that the above inversion formulae do not require the model to be identifiable. We saw in Section 2 that a crucial factor in obtaining identifiability is using the symmetry of the unknown mixed distribution. In the next paragraph we use the symmetry of the mixed distribution to provide a contrast function.

3.2. *A contrast function.* The second key point follows from the following simple remark. Let $F_\theta = A_\theta^{-1} G = A_\theta^{-1} A_{\theta_0} F_0$, where $\theta \in \Theta$. Clearly, if $\theta = \theta_0$, we have $F_\theta = F_0$ (from Section 3.1), and it must have the invariance property of c.d.f.'s of symmetric distributions, $F_0(x) = 1 - F_0(-x)$. For simplicity, let us introduce $S_r$, the *symmetry operator* defined by $S_r F(\cdot) = 1 - F(-\cdot)$. The preceding remark may be reformulated as follows: if $\theta = \theta_0$, then $A_\theta^{-1} G = S_r A_\theta^{-1} G$ or, equivalently, $G = A_\theta S_r A_\theta^{-1} G$, by applying $A_\theta$ on the left-hand side of the last equality. What about the converse? The answer is given in the following theorem, whose proof is given in Section 5.

THEOREM 3.1. *Consider model* (2) *with $F_0$ the c.d.f. of a symmetric distribution and $\theta_0 \in \Theta$. If, for $\theta \in \Theta$, we have $G = A_\theta S_r A_\theta^{-1} G$, then $\theta = \theta_0$.*

Assuming that $G$ is known, we can recover the true value $\theta_0$ of $\theta$ by minimizing a discrepancy measure between $G$ and $G_\theta = A_\theta S_r A_\theta^{-1} G$. Recall that $G$ is unknown but can be estimated, which is why we choose to consider the discrepancy measure $K$, defined by

$$(11) \qquad K(\theta) \equiv K(\theta; G) = \int_{\mathbb{R}} (G_\theta(x) - G(x))^2 \, dG(x), \qquad \theta \in \Theta.$$

The choice of introducing the weighted measure $G$ in the above integral follows from the consideration that if $G$ is replaced by its empirical c.d.f., then the integral sign turns into a simple sum. As a consequence of the preceding theorem, assuming that $F$ is sufficiently smooth and that $G$ is known, we are able to show that $K$ is a contrast function for the unknown Euclidean parameter $\theta$.

COROLLARY 3.1. *Under assumption* C1, *$K$ is a contrast function: for all $\theta \in \Theta$, $K(\theta) \geq 0$ and $K(\theta) = 0$ if and only if $\theta = \theta_0$.*

3.3. *Estimators of the Euclidean parameter $\theta$.* The above Corollary 3.1 suggests that the unknown Euclidean parameter $\theta$ should be estimated as follows:

$$\hat{\theta}_n = \underset{\theta \in \Theta}{\arg\min}\, K(\theta; \hat{G}_n),$$

where $\hat{G}_n$ is an estimator of the c.d.f. of $G$. It is important to note that if $\hat{G}_n$ is a stepwise function, $K(\theta; \hat{G}_n)$ is also a stepwise function with respect to parameters $\mu_1$ and $\mu_2$, and does not have the required regularity properties for differentiable optimization techniques to be applied in order to find $\hat{\theta}_n$. This is the reason why we need to distinguish two cases: (P1) the parameters $\mu_1$ and $\mu_2$ are known and (P2) the parameters $\mu_1$ and $\mu_2$ are unknown.

(P1) The parameters $\mu_1$ and $\mu_2$ are known, whereas $\lambda$ and $F$ are unknown.

For this problem, we suppose that the true mixing proportion $\lambda_0$ belongs to $[0, 1/2 - d]$, where $d \in (0, 1/2)$. In this case the parameter $\theta$ reduces to $\lambda$ and we estimate $\lambda$ by

$$\hat{\lambda}_n = \underset{\lambda \in [0, 1/2 - d]}{\arg\min}\, K(\lambda; \hat{G}_n),$$

where $\hat{G}_n$ is the empirical c.d.f. of $G$ defined by

$$(12) \qquad \hat{G}_n(x) = \frac{1}{n} \sum_{i=1}^{n} \mathbf{1}_{X_i \leq x} \qquad \forall\, x \in \mathbb{R},$$

where $\mathbf{1}$ denotes the indicator function. Let us give an explicit formula for $\hat{G}_\lambda^{(n)} = A_\lambda S_r A_\lambda^{-1} \hat{G}_n$ involving a sum of $n$ terms:

$$(13) \quad \hat{G}_\lambda(x) = 1 + \frac{1}{n} \sum_{i=1}^{n} \left( \frac{\lambda}{\lambda - 1} \mathbf{1}_{x \leq \eta + 2\mu_1 - X_i} + \frac{1 - 2\lambda}{\lambda} \left( \frac{\lambda}{\lambda - 1} \right)^{L(i,x)} \right),$$

where

$$L(i, x) = \max \left( 1, \left\lceil \frac{x - 2\mu_1 + X_i - \eta}{\eta} \right\rceil \right),$$

and where $\lceil x \rceil$ denotes the smallest integer greater than or equal to $x$ and $\eta = \mu_2 - \mu_1$. The following theorem, whose proof is provided in Section 5, gives the asymptotic behavior of $\hat{\lambda}_n$.

THEOREM 3.2. *Assume that the c.d.f. $F_0$ satisfies assumption* C1. *Then* (i) $\hat{\lambda}_n$ *converges almost surely to $\lambda_0$, and* (ii) *we have $\sqrt{n}(\hat{\lambda}_n - \lambda_0) = O_P(1)$.*

Note that if $F_0$ is assumed to admit a first-order moment, then, using the first-order moment equation of $g$, we show that $\lambda_0$ can be directly estimated by the natural empirical estimator

$$\bar{\lambda}_n = \frac{n^{-1}\sum_{i=1}^n X_i - \mu_2}{\mu_2 - \mu_1},$$

which obviously satisfies results of the above theorem.

(P2) The parameters $\mu_1, \mu_2, \lambda$ and $F$ are unknown.

For this problem, we suppose that $\Theta = [0, 1/2 - d] \times \mathcal{X}$, where $0 < d < 1/2$ and $\mathcal{X}$ is a compact subset of $\mathbb{R}^2$ such that $\mathcal{X} \cap \Delta = \varnothing$, and the unknown Euclidean parameter $\theta$ is an interior point of $\Theta$. As explained previously, we need to change $K(\cdot; \hat{G}_n)$ into the more regular version $K_r(\cdot; \hat{G}_n)$ defined by

$$K_r(\theta; \hat{G}_n) = \int_{\mathbb{R}} (\tilde{G}_\theta^{(n)}(x) - \tilde{G}_n(x))^2 \, d\hat{G}_n(x),$$

where $\tilde{G}_\theta^{(n)} = A_\theta S_r A_\theta^{-1} \tilde{G}_n$ and $\tilde{G}_n(x) = \int_{-\infty}^x \hat{g}_n(y) \, dy$, with

$$\hat{g}_n(x) = \frac{1}{b_n} \int_{\mathbb{R}} q\left(\frac{x-y}{b_n}\right) d\hat{G}_n(y),$$

where $(b_n)_{n \geq 1}$ is a sequence of real numbers decreasing to 0. Our numerical applications are based upon the kernel function $q$ defined by $q(x) = (1 - |x|)\mathbf{1}_{|x| \leq 1}$.

As for the (P1) problem, we prove in Section 5 asymptotic results summarized in the next theorem for the estimator $\hat{\theta}_n$. From a general point of view, $\tilde{G}_n$ is a smooth estimate of the c.d.f. $G$ defined, for $x \in \mathbb{R}$, by

$$\tag{14} \tilde{G}_n(x) = \frac{1}{n}\sum_{k=1}^n Q\left(\frac{x - X_k}{b_n}\right),$$

where $Q(x) = \int_{-\infty}^x q(y) \, dy$, with $q$ an even density function with compact support and second-order moment equal to 1, and $(b_n)_{n \geq 1}$ is a sequence of nonnegative real numbers decreasing to 0 with $nb_n \to +\infty$ and $\sqrt{n}b_n^2 = O(1)$. The fact that $q$ has compact support leads to an explicit formula for $\tilde{G}_\theta^{(n)}$, involving a sum of $n$ terms,

$$\tilde{G}_\theta(x) = 1 + \frac{1}{(n\,b_n)}\sum_{i=1}^n \left\{ \frac{\lambda}{\lambda - 1} Q\left(\frac{-x + \eta + 2\mu_1 - X_i}{b_n}\right) \right.$$
$$\left. + \frac{1 - 2\lambda}{\lambda}\left(\frac{\lambda}{\lambda - 1}\right)^{L_2(i,x)} \right.$$
$$\tag{15}$$

$$+ \frac{2\lambda - 1}{\lambda(\lambda - 1)} \sum_{k=L_1(i,x)}^{L_2(i,x)-1} \left( \frac{\lambda}{\lambda - 1} \right)^k$$

$$\times Q\left( \frac{-x + (k+1)\eta + 2\mu_1 - X_i}{b_n} \right) \Bigg\},$$

where, for $k = 1, 2$,

$$L_k(i,x) = \max\left( 1, \left\lceil \frac{x - 2\mu_1 + X_i - \eta + (-1)^k b_n}{\eta} \right\rceil \right).$$

THEOREM 3.3. *If the c.d.f. $F_0$ satisfies* C1, *then $\hat{\theta}_n$ converges almost surely to $\theta_0$. If, in addition, $F_0$ satisfies* C2, *we have $n^{1/4-\alpha}(\hat{\theta}_n - \theta_0) = o_{\text{a.s.}}(1)$, for all $\alpha > 0$.*

3.4. *Estimators of functional parameter $F$.* As suggested by the inversion formula (8), once we get a consistent estimator $\hat{\theta}_n$ of the unknown (or partially unknown) Euclidean true parameter $\theta_0$, it is natural to seek to approximate the unknown c.d.f. $F_0$ by $\tilde{F}_n = A_{\hat{\theta}_n}^{-1} \hat{G}_n$. However, since we approximate the c.d.f. of a symmetric distribution, we constrain $\tilde{F}_n$ to satisfy the invariance property $\tilde{F}_n = S_r \tilde{F}_n$, leading to the final estimator

$$(16) \qquad \hat{F}_n = \tfrac{1}{2}(I + S_r) A_{\hat{\theta}_n}^{-1} \hat{G}_n,$$

where $I$ is the identity operator. By similar arguments, the unknown density function $f_0$ can in turn be estimated by

$$(17) \qquad \hat{f}_n = \tfrac{1}{2}(I + S_d) A_{\hat{\theta}_n}^{-1} \hat{g}_n,$$

where the operator $S_d$ is defined by $(S_d f)(x) = f(-x)$ (corresponding to the invariance property of densities of symmetric distributions). The next theorem gives asymptotic results for both $\hat{F}_n$ and $\hat{f}_n$ for problems (P1) and (P2). These theorems are proved in Section 5.

THEOREM 3.4. (i) *If $F_0$ satisfies* C1, *then $\|\hat{F}_n - F_0\|_\infty$ converges almost surely to 0 for problems* (P1) *and* (P2).

(ii) *Under* C1, *we have $\|\hat{F}_n - F_0\|_\infty = O_P(n^{-1/2})$ for problem* (P1). *Under* C2, *for problem* (P2) *we have $\|\hat{F}_n - F_0\|_\infty = o_{\text{a.s.}}(n^{-1/4+\alpha})$ for any $\alpha > 0$.*

(iii) *Under* C1 (*resp.* C2*), for problem* (P1) [*resp.* (P2)], *$\|\hat{f}_n - f_0\|_1$ converges almost surely to 0.*

Let us notice that generally $\hat{F}_n$ (resp. $\hat{f}_n$) is not a c.d.f. function (resp. a density function). Indeed, by the definition of a mixture, $g$ belongs to the range of the operator $A_\theta$, whereas this is no longer true for its approximate $\hat{g}_n$. Since there is no possibility that $\hat{g}_n$ is a two-component mixture in the sense of model (2), it follows that $A_{\hat{\theta}_n}^{-1} \hat{g}_n$ cannot be a density function, and the same holds for $\hat{f}_n$. However, from a practical point of view, we can easily transform estimators $\hat{f}_n$ into density functions. Let us consider $f_n^* = \hat{f}_n \mathbf{1}_{\hat{f}_n \geq 0}$. It is straightforward to show that $\|f_n^* - f_0\|_1 \leq \|\hat{f}_n - f_0\|_1$, and then we have the almost sure convergence of $\|f_n^* - f_0\|_1$ to 0, given the assumptions of Theorem 3.4(iii). Moreover, under the same assumptions and with $s_n = \int_{\mathbb{R}} f_n^*(x)\, dx$, we have

$$|s_n - 1| = \left| \int_{\mathbb{R}} (f_n^*(x) - f_0(x))\, dx \right|$$
$$\leq \|f_n^* - f_0\|_1$$
$$\leq \|\hat{f}_n - f_0\|_1 \to 0 \qquad \text{a.s.}$$

Therefore, $\tilde{f}_n = s_n^{-1} f_n^*$ are density functions that satisfy $\|\tilde{f}_n - f_0\|_1 \to 0$, almost surely.

3.5. *Discussion of the three-component case.* As we discussed in Section 2, identifiability results exist (see [19]) for the following three-component model:

$$G(x) = \lambda_1 F(x - \mu_1) + \lambda_2 F(x - \mu_2) + \lambda_3 F(x - \mu_3) \qquad \forall x \in \mathbb{R},$$

where $F$ is the c.d.f. of a symmetric distribution, and the $\lambda_i$'s are nonnegative real numbers with $\lambda_1 + \lambda_2 + \lambda_3 = 1$. A question naturally arises concerning the possibility of extending our estimation method to the above model. Following the method presented in Section 3.1, we get, for all $\ell \geq 1$,

$$F(x) = \frac{G(x + \mu_3)}{\lambda_3}$$
$$+ \sum_{k=1}^{\ell-1} (-1)^k \sum_{(i_1, \ldots, i_k) \in \{1,2\}^k} \frac{\lambda_{i_1} \cdots \lambda_{i_k}}{\lambda_3^{k+1}} G(x + \mu_3 + \eta_{i_1} + \cdots + \eta_{i_k})$$
$$+ (-1)^\ell \sum_{(i_1, \ldots, i_\ell) \in \{1,2\}^\ell} \frac{\lambda_{i_1} \cdots \lambda_{i_\ell}}{\lambda_3^\ell} F(x + \eta_{i_1} + \cdots + \eta_{i_\ell}) \qquad \forall x \in \mathbb{R},$$

where we suppose that $\lambda_3 > \max(\lambda_1, \lambda_2)$ and we denote $\eta_i = \mu_3 - \mu_i$ for $i = 1, 2$. To prove that a type (8) formula exists, we need to show that, for all $x \in \mathbb{R}$, we have

$$
\begin{aligned}
(18) \quad F(x) = {} & \frac{G(x + \mu_3)}{\lambda_3} \\
& + \sum_{k=1}^{+\infty} (-1)^k \sum_{(i_1, \ldots, i_k) \in \{1,2\}^k} \frac{\lambda_{i_1} \cdots \lambda_{i_k}}{\lambda_3^{k+1}} G(x + \mu_3 + \eta_{i_1} + \cdots + \eta_{i_k}).
\end{aligned}
$$

Unfortunately, taking $x \geq 1$, it is easy to see that (18) is not satisfied by taking, for example, $\lambda_1 = \lambda_2 = 4/15$, $\lambda_3 = 7/15$, $\mu_1 = 0$, $\mu_2 = -1$, $\mu_3 = 1$ and $F$ c.d.f. of the uniform distribution on $(-1, 1)$. Note, however, that if the inversion formula (18) is valid [this is the case, e.g., for $2 \max(\lambda_1, \lambda_2) < \lambda_3$], the methodology proposed in this section for the two-component case may be applied.

**4. Numerical study.** We consider two distinct problems. The first is to estimate $\lambda$ given that $\mu_1$ and $\mu_2$ are known. In this case we use an explicit formula for $\hat{G}_\lambda^{(n)}$. In the second case we estimate $\theta = (\lambda, \mu_1, \mu_2)$ and we consider $\tilde{G}_\theta^{(n)}$, the regularized version of $\hat{G}_\theta^{(n)}$. Explicit formulae for $\hat{G}_\lambda^{(n)}$ and $\tilde{G}_\theta^{(n)}$ are given in (13) and (15). Recall that the computation of $\tilde{G}_\theta^{(n)}$ involves the choice of a bandwidth $b_n$. All the simulation results have been obtained with $b_n = n^{-1/4}$. This value is not optimal to estimate the density $g$ but it is compatible with the assumption $\sqrt{n} b_n^2 = O(1)$ needed to achieve the convergence rate given in Theorem 3.3. Note that in all our simulations the variance $\sigma_g^2$ under $g$ is close to 1; our choice for $b_n$ is then close to the bandwidth that minimizes the mean integrated squared error, usually approximated by $\sigma_g(4/3n)^{1/5}$ (see, e.g., [3]). It is known to be a good approximation for normal data and a Gaussian kernel but we cannot insure that it leads to an optimal choice for our problem. For the real example of rainfall data given at the end of this section, we used the bandwidth ($b_n = 3.84$) provided by the R software.

*Choice of optimization method.* Problem (P1) attempts to find an estimate $\hat{\lambda}_n$ of $\lambda$ when $\mu_1$ and $\mu_2$ are known,

$$
(19) \qquad \hat{\lambda}_n = \underset{\lambda \in [0, 1/2 - d]}{\arg\min} \, K(\lambda; \hat{G}_n).
$$

Problem (P2) attempts to find an estimate $\hat{\theta}_n$ of $\theta = (\lambda, \mu_1, \mu_2)$,

$$
(20) \qquad \hat{\theta}_n = \underset{\theta \in \Theta}{\arg\min} \, \tilde{K}_r(\theta; \hat{G}_n).
$$

Both problems require the minimization of a differentiable functional. As far as problem (19) is concerned, numerical experiments indicate that $K(\cdot; \hat{G}_n)$ is strictly convex in $[0, 1/2 - d]$ and, thus, an unconstrained minimization algorithm can safely be used, with a starting point in this interval. We use the quasi-Newton BFGS (Broyden, Fletcher, Goldfarb and Shanno) method (see, e.g., [32]). In the second case, some experiments with the same unconstrained method show that $K_r(\cdot; \hat{G}_n)$ is not convex, and that $K_r(\cdot; \hat{G}_n)$ has local minima not belonging to $\Theta$. So we use the constrained version of the BFGS algorithm, where bounds on the variables can be taken into account. In both cases, we provide the gradient of the functional, which can be readily computed from the explicit formulae given in Section 3.3. All the computations are performed with Scilab.

TABLE 1
*Empirical mean and standard error (in brackets) of $\lambda$ estimates, obtained from 500 simulations of i.i.d. samples of size n, for the (P1) problem with $\mu_1 = -1$ and $\mu_2 = 2$*

| $n$ $\lambda$ | 0.15 | 0.25 | 0.35 |
|---|---|---|---|
| 100 | 0.151 (0.058) | 0.256 (0.060) | 0.347 (0.057) |
| 400 | 0.148 (0.031) | 0.252 (0.032) | 0.349 (0.029) |

*Numerical result of the Monte Carlo study for Gaussian mixtures.* In this section we denote by $\mathcal{N}(\mu, \sigma^2)$ a Gaussian distribution with mean $\mu$ and variance $\sigma^2$. The performance of our method is evaluated, via a Monte Carlo study, on the Gaussian mixture

$$\lambda * \mathcal{N}(\mu_1, 1) + (1 - \lambda) * \mathcal{N}(\mu_2, 1),$$

for the (P1) problem (see Table 1) and in the (P2) problem (see Tables 2 and 3). More precisely, Table 1 summarizes the performance of our method for different values of $\lambda$, that is, $\lambda = 0.15$ (weakly bumped model), $\lambda = 0.25$ (moderately bumped model) and $\lambda = 0.35$ (strongly bumped model), when

TABLE 2
*Empirical mean and standard error of $(\lambda, \mu_1, \mu_2)$ semiparametric estimates, obtained from 200 simulations of i.i.d. samples of size n, for the (P2) problem with $b_n = n^{-1/4}$*

| $n$ | $(\lambda, \mu_1, \mu_2)$ | Empirical means | Standard errors |
|---|---|---|---|
| 100 | $(0.15, -1, 2)$ | $(0.161, -0.948, 2.030)$ | $(0.052, 0.365, 0.137)$ |
| 200 | $(0.15, -1, 2)$ | $(0.157, -1.027, 2.023)$ | $(0.035, 0.283, 0.101)$ |
| 100 | $(0.25, -1, 2)$ | $(0.249, -1.011, 2.009)$ | $(0.060, 0.289, 0.154)$ |
| 200 | $(0.25, -1, 2)$ | $(0.251, -1.000, 2.010)$ | $(0.041, 0.195, 0.101)$ |
| 100 | $(0.35, -1, 2)$ | $(0.347, -0.988, 1.990)$ | $(0.056, 0.230, 0.145)$ |
| 200 | $(0.35, -1, 2)$ | $(0.357, -0.976, 2.012)$ | $(0.046, 0.176, 0.114)$ |

Table 3

*Empirical mean and standard error of $(\lambda, \mu_1, \mu_2)$ maximum likelihood estimates, obtained from 200 simulations of i.i.d. samples of size $n$, for the (P2) problem with $b_n = n^{-1/4}$*

| $n$ | $(\lambda, \mu_1, \mu_2)$ | **Empirical means** | **Standard errors** |
|-----|------|------|------|
| 100 | $(0.15, -1, 2)$ | $(0.163, -0.987, 2.018)$ | $(0.054, 0.431, 0.138)$ |
| 200 | $(0.15, -1, 2)$ | $(0.152, -1.013, 2.004)$ | $(0.035, 0.283, 0.089)$ |
| 100 | $(0.25, -1, 2)$ | $(0.256, -1.008, 2.020)$ | $(0.051, 0.268, 0.132)$ |
| 200 | $(0.25, -1, 2)$ | $(0.247, -1.003, 2.004)$ | $(0.046, 0.204, 0.114)$ |
| 100 | $(0.35, -1, 2)$ | $(0.342, -1.041, 1.980)$ | $(0.054, 0.231, 0.161)$ |
| 200 | $(0.35, -1, 2)$ | $(0.345, -1.009, 1.991)$ | $(0.041, 0.159, 0.111)$ |

$\mu_1 = -1$ and $\mu_2 = 2$ are known. Table 2 summarizes the performance of our method in estimating $\lambda = 0.15, 0.25, 0.35$, and $\mu_1 = -1$, $\mu_2 = 2$, while Table 3 gives the performance of the standard maximum likelihood approach in the same framework.

*Comments on Tables* 1–3. The results in Table 1 show first that empirical bias amounts to less than 1% of the true values, and that standard errors are reasonably small. In order to clarify the analysis of the results given in Table 1 and to quantify the influence of bumps on the estimation efficiency, we can normalize the empirical standard errors with respect to the true values of the parameters (std/$\lambda$). We obtain for $\lambda = 0.15, 0.25, 0.35$, normalized empirical standard errors equal to 0.386, 0.240, 0.162, respectively, for $n = 100$, and equal to 0.206, 0.128, 0.0828, for $n = 400$. These indicators show, roughly speaking, that our estimation method is around 2.4 times more precise when $\lambda = 0.35$, in comparison with the case where $\lambda = 0.15$, and 1.6 times more precise when $\lambda = 0.35$ in comparison with the case where $\lambda = 0.25$. This shows that the nonnegligibility of one subpopulation with respect to the other subpopulation improves the quality of the estimators.

Concerning Tables 2 and 3, it is interesting to note that, when the location parameters are unknown, the previous remark is no longer true. In fact, even if the previous comments on empirical standard errors are clearly relevant, it is worth noting that the smaller empirical bias is not obtained for the highly bumped model, but for the moderately bumped model. To explain this phenomenon, we can remark that when $\lambda = 0.15$, there are few data to estimate $\mu_1$, whereas when $\lambda = 0.35$, even if there are many more data to estimate $\mu_1$, this estimation is disturbed by the left tail of the distribution centered on $\mu_2$. Finally, it is with $\lambda = 0.25$ that we obtain the best compromise and therefore the best estimates with regard to minimum bias. In addition, we observe that the performance of the maximum likelihood approach (which is known to be asymptotically efficient) is in the range of those obtained by our method, which illustrates the good behavior of our semiparametric approach with respect to the parametric approach.

*A trimodal example.* We use a basic symmetric density $f$ which is already a mixture, that is,

$$f(x) = \tfrac{1}{8}\varphi(x+4) + \tfrac{3}{4}\varphi(x) + \tfrac{1}{8}\varphi(x-4),$$

where $\varphi$ is the density function of the standard Gaussian distribution. The density of the simulated data is taken as

$$g(x) = \tfrac{1}{4}f(x) + \tfrac{3}{4}f(x-4) \qquad \forall\, x \in \mathbb{R}.$$

We performed the estimation on a simulated sample of size $n = 100$. The results are given in Figure 1. Figure 1(a) shows $\tilde{f}$ superimposed with the true density function $f$ and Figure 1(b) shows the reconstructed density function $\tilde{g}(\cdot) = \hat{\lambda}\tilde{f}(\cdot - \hat{\mu}_1) + (1-\hat{\lambda})\tilde{f}(\cdot - \hat{\mu}_2)$, using estimated values of $\lambda$, $\mu_1$ and $\mu_2$, superimposed with the true density function $g$. The optimization required 31 iterations and 45 evaluations of $K_r(\cdot; \hat{G}_n)$ and its gradient.

Standard errors for Euclidean parameters are computed by the Jackknife method (see, e.g., [14]). We observe that for a reasonably small sample size $n = 100$ the reconstructed mixture density $\tilde{g}$ almost yields the true density $g$. The main differences appear around local modes and in the tails of $g$.

*Numerical results on real data.* We use the average amount of precipitation (rainfall) in inches for each of 70 United States (and Puerto Rican) cities (from the Statistical Abstract of the United States, 1975; see [30]). We consider two models. The first is model (2) in which we denote by $\hat{\lambda}$, $\hat{\mu}_1$, $\hat{\mu}_2$ and $\tilde{f}$ the estimators of $\lambda$, $\mu_1$, $\mu_2$ and $f$ (the density function of the c.d.f. $F$). The second model is a parametric version of model (2) in which we assume that $f$ is the density function of a centered Gaussian distribution with variance equal to $\sigma^2$. Estimators of unknown parameters of the second model are denoted by $\tilde{\lambda}$, $\tilde{\mu}_1$, $\tilde{\mu}_2$ and $\tilde{\sigma}^2$, and calculated according to the maximum likelihood method.

Figure 2(a) shows $\tilde{f}$, the estimator of $f$ superimposed with the density function of $\mathcal{N}(0, \tilde{\sigma}^2)$. Figure 2(b) shows $\hat{g}$, the empirical estimate of $g$ (obtained by the kernel method) superimposed with both the reconstructed density $\hat{\lambda}\tilde{f}(\cdot - \hat{\mu}_1) + (1-\hat{\lambda})\tilde{f}(\cdot - \hat{\mu}_2)$ (using estimated values $\hat{\lambda}$, $\hat{\mu}_1$ and $\hat{\mu}_2$ of $\lambda$, $\mu_1$ and $\mu_2$) and the density of the parametric model where the Euclidean parameter is replaced by its maximum likelihood estimator. The optimization required 32 iterations and 66 evaluations of $K_r(\cdot; \hat{G}_n)$ and its gradient.

We observe in Figure 2(a) that the nonparametric density estimate $\hat{f}$ is provided with two symmetric small bumps at the beginnings of its tails, while the best fitting Gaussian density does not obviously benefit from this kind of singularity and is sharper around the origin. In Figure 2(b) we can see that these two additional bumps make the difference in the good fitting behavior of the reconstructed mixing distribution, except in a small area around $[-20, 0]$ (the area of interest being $[-20, 75]$), where the best fitting Gaussian mixture is slightly closer to $\hat{g}$. Notice also that the smaller bump on the left of $\hat{g}$ is clearly detected by our method, while the best fitting Gaussian mixture almost misses this singularity. Again, standard errors (given in brackets) for Euclidean parameters are computed using the Jackknife method.

## 5. Proofs.

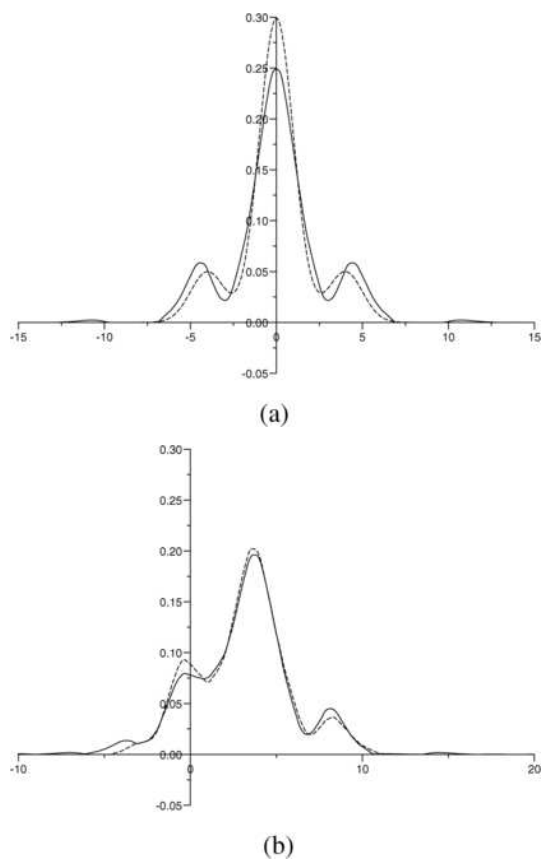

(a)

(b)

Fig. 1. *Estimated parameters* $\hat{\mu}_1 = 0.691$ *(0.760),* $\hat{\mu}_2 = 3.728$ *(0.153),* $\hat{\lambda} = 0.232$ *(0.079) for* $n = 100$ *and* $b_n = n^{-1/4}$, *the results in parentheses corresponding to the empirical standard errors.* (a) *Graph of* $\hat{\tilde{f}}$ *(solid) and graph of* $f$ *(dashed).* (b) *Graph of* $\hat{\lambda}\hat{f}(\cdot - \hat{\mu}_1) + (1 - \hat{\lambda})\hat{f}(\cdot - \hat{\mu}_2)$ *(solid) and graph of* $g$ *(dashed).*

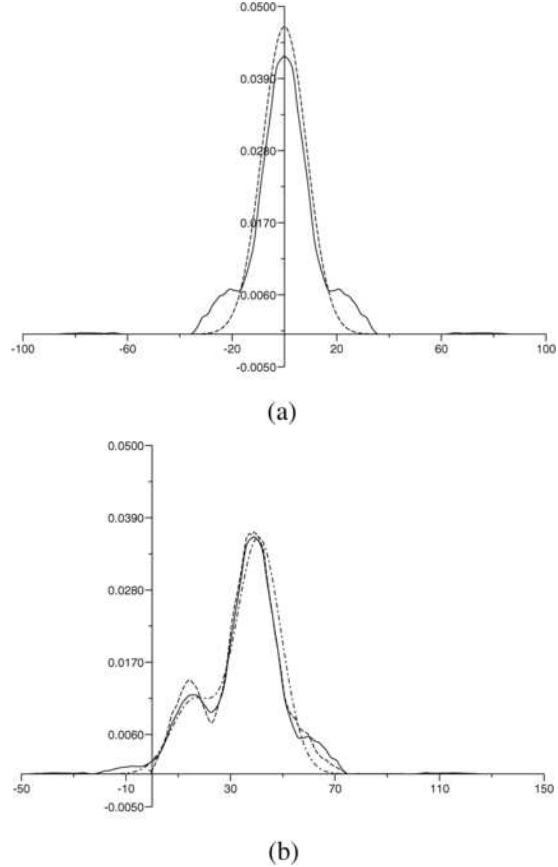

**(a)**

**(b)**

Fig. 2. *Estimated parameters for model* (2): $\hat{\mu}_1 = 13.107$ (3.299), $\hat{\mu}_2 = 39.056$ (1.395), $\hat{\lambda} = 0.171$ (0.078) *(the bandwidth is fixed at 3.84). Estimated parameters for model* $\lambda * \mathcal{N}(\mu_1, \sigma^2) + (1 - \lambda) * \mathcal{N}(\mu_2, \sigma^2)$: $\tilde{\mu}_1 = 15.715$ (2.220), $\tilde{\mu}_2 = 40.773$ (1.297), $\tilde{\lambda} = 0.235$ (0.060) *and* $\tilde{\sigma} = 8.504$ (1.187), *the results in parentheses corresponding to the empirical standard errors.* (a) *Graph of the nonparametric density estimator* $\hat{f}$ *and graph of the density of* $\mathcal{N}(0, \tilde{\sigma}^2)$ *(dashed).* (b) *Graph of* $\hat{g}$ *(dashed), graph of* $\tilde{\lambda}\hat{f}(\cdot - \hat{\mu}_1) + (1 - \tilde{\lambda})\hat{f}(\cdot - \hat{\mu}_2)$ *(solid) and graph of the density of* $\tilde{\lambda}\mathcal{N}(\tilde{\mu}_1, \tilde{\sigma}^2) + (1 - \tilde{\lambda})\mathcal{N}(\tilde{\mu}_2, \tilde{\sigma}^2)$ *(dash-dot).*

5.1. *Notation and preliminary results.* According to whether we are looking at density function estimation or c.d.f. estimation, the operators $A_\theta$ and $A_\theta^{-1}$, given in (10), are defined, respectively, on spaces $L^1(\mathbb{R})$ or $L^\infty(\mathbb{R})$ (endowed with the usual norms $\|\cdot\|_1$ and $\|\cdot\|_\infty$, resp.). Independently of the space under consideration, it is straightforward to check that the norms (denoted $\|\cdot\|$) of operators $A_\theta$ and $A_\theta^{-1}$, for $\lambda \in [0, 1/2 - d]$ and $d \in (0, 1/2)$,

satisfy

$$(21) \qquad \|A_\theta\| \le 1 \quad \text{and} \quad \|A_\theta^{-1}\| \le \frac{1}{1-2\lambda} \le \frac{1}{2d}.$$

Let us recall some basic results on $\hat{G}_n$ and $\tilde{G}_n$ defined respectively by (12) and (14). From well-known results on empirical processes (see, e.g., [37]), for general distribution functions $G$, we have

$$(22) \qquad \sqrt{n}\|\hat{G}_n - G\|_\infty = O_P(1),$$

and the law of iterated logarithm (LIL)

$$(23) \qquad \|\hat{G}_n - G\|_\infty = O_{\text{a.s.}}\left( \sqrt{\frac{\log\log(n)}{n}} \right).$$

If $\|f\|_\infty < \infty$, and if $f$ has derivative $f^{(1)}$ with $\|f^{(1)}\|_\infty < \infty$, the same holds for $g$, and by Corollary 1, page 766 in [37], if $q$ has compact support, and if $\sqrt{n}b_n^2 = O(1)$, then we have

$$(24) \qquad \sqrt{n}\|\hat{G}_n - \tilde{G}_n\|_\infty = O_{\text{a.s.}}(1).$$

Hence, the result (23) holds for $\tilde{G}_n$.

In the remainder of this paper we denote by $\dot{L}$ and $\ddot{L}$ the first- and second-order derivatives of a general function $L$ with respect to $\lambda \in [0, 1/2 - d]$ for problem (P1) and $\theta \in \Theta = [0, 1/2 - d] \times \mathcal{X}$ for problem (P2) (see Section 3.3 for assumptions on the Euclidean parameter space). In the sequel $|\cdot|_2$ denotes the Euclidean norm.

LEMMA 5.1. *There exists $c \in (0, +\infty)$, such that for all $\theta \in \Theta$ and $n \ge 1$, we have*

$$(25) \qquad \|\hat{G}_\theta^{(n)} - G_\theta\|_\infty \le c\|\hat{G}_n - G\|_\infty.$$

PROOF. Straightforward, since we have

$$\begin{aligned}
\|\hat{G}_\theta^{(n)} - G_\theta\|_\infty &= \|A_\theta S_r[A_\theta^{-1}(\hat{G}_n - G)]\|_\infty \\
&\le \|A_\theta^{-1}\| \times \|\hat{G}_n - G\|_\infty \\
&\le \frac{1}{2d}\|\hat{G}_n - G\|_\infty. \qquad \square
\end{aligned}$$

LEMMA 5.2. *Under C1, the mapping $\theta \mapsto G_\theta(x)$ is Lipschitz on $\Theta$ uniformly in $x \in \mathbb{R}$, and the contrast function $K$ is Lipschitz on $\Theta$.*

PROOF.   For all $(\theta, \theta') \in \Theta^2$, we have $|K(\theta) - K(\theta')| \leq C \|G_\theta - G_{\theta'}\|_\infty$. Therefore, it is sufficient to prove that $\theta \mapsto G_\theta(\cdot)$ is uniformly Lipschitz on $\Theta$. Simple calculations lead to

$$(26) \quad \|G_\theta - G_{\theta'}\|_\infty \leq \|A_\theta S_r A_\theta^{-1} G - A_{\theta'} S_r A_{\theta'}^{-1} G\|_\infty + \|A_\theta^{-1} G - A_{\theta'}^{-1} G\|_\infty.$$

First we prove that the first term on the right-hand side of (26) is Lipschitz. Let us remark now that, for all bounded functions $H$, we have, for all $x \in \mathbb{R}$,

$$(27) \quad \begin{aligned} |A_\theta H(x) - A_{\theta'} H(x)| &\leq 2|\lambda - \lambda'| \times \|H\|_\infty \\ &\quad + \sup_{x \in \mathbb{R}} \max_{i=1,2} |H(x) - H(x - (\mu_i - \mu_i'))|. \end{aligned}$$

On the other hand, noticing $\eta = \mu_2 - \mu_1$ (resp. $\eta' = \mu_2' - \mu_1'$), we remark that for all $\theta \in \Theta$, $A_\theta^{-1} G$ satisfies, for all $(z, z') \in \mathbb{R}^2$,

$$\begin{aligned} &|S_r A_\theta^{-1} G(z) - S_r A_\theta^{-1} G(z')| \\ &\quad = |A_\theta^{-1} G(-z) - A_\theta^{-1} G(-z')| \\ &\quad = \left| \frac{1}{1-\lambda} \sum_{k=0}^\infty \left( -\frac{\lambda}{1-\lambda} \right)^k (G(-z + \mu_2 + k\eta) - G(-z' + \mu_2 + k\eta)) \right| \\ &\quad \leq \frac{1}{2d} \sup_{y \in \mathbb{R}} |G(y) - G(y - z + z')| \\ &\quad \leq \frac{1}{2d} |G|_{\mathrm{Lip}} |z - z'|, \end{aligned}$$

because under C1 $G$ is Lipschitz with Lipschitz constant $|G|_{\mathrm{Lip}}$. Now replacing $H$ in (27) by $S_r A_\theta^{-1} G$, it follows from the above inequality that

$$(28) \quad |A_\theta S_r A_\theta^{-1} G(x) - A_{\theta'} S_r A_{\theta'}^{-1} G(x)| \leq C|\theta - \theta'|_2.$$

It remains to be proved that the second term on the right-hand side of inequality (26) is Lipschitz. We have

$$\begin{aligned} &|A_\theta^{-1} G(x) - A_{\theta'}^{-1} G(x)| \\ &\quad \leq \left| \frac{1}{1-\lambda} \sum_{k=0}^\infty \left( -\frac{\lambda}{1-\lambda} \right)^k (G(x + \mu_2 + k\eta) - G(x + \mu_2' + k\eta')) \right| \\ &\quad\quad + \left| \left( \frac{1}{1-\lambda} \sum_{k=0}^\infty \left( -\frac{\lambda}{1-\lambda} \right)^k - \frac{1}{1-\lambda'} \sum_{k=0}^\infty \left( -\frac{\lambda'}{1-\lambda'} \right)^k \right) \right. \\ &\quad\quad\quad\quad \left. \times G(x + \mu_2 + k\eta) \right|. \end{aligned}$$

$G$ is supposed to be Lipschitz. We have, for all $x \in \mathbb{R}$,

$$|G(x + \mu_2 + k\eta) - G(x + \mu_2' + k\eta')| \leq |G|_{\mathrm{Lip}} (k+1)|\theta - \theta'|_2;$$

thus, we obtain, by the two previous inequalities,

$$
\begin{aligned}
&|A_\theta^{-1}G(x) - A_{\theta'}^{-1}G(x)| \\
(29) \qquad &\leq \frac{c_1}{1-\lambda'}\sum_{k=0}^{\infty}(k+1)\left(\frac{\lambda'}{1-\lambda'}\right)^k|\theta-\theta'|_2 + c_2\|G\|_\infty|\lambda-\lambda'| \\
&\leq c_3|\theta-\theta'|_2,
\end{aligned}
$$

where $c_1$, $c_2$ and $c_3$ are nonnegative real constants. From inequalities (26)–(29), it follows that the function $\theta \mapsto G_\theta(x)$ is Lipschitz on $\Theta$ uniformly in $x \in \mathbb{R}$, and thus, $K$ is Lipschitz on $\Theta$. $\quad\square$

LEMMA 5.3. *For any $\alpha > 0$, under* C1 *we have*

$$
(30) \qquad \sup_{\theta\in\Theta}|K(\theta;\hat{G}_n) - K(\theta)| = o_{\mathrm{a.s.}}(n^{-1/2+\alpha}).
$$

The same result holds replacing $K(\cdot;\hat{G}_n)$ by $K_r(\cdot;\hat{G}_n)$. It is an obvious consequence of properties of $\tilde{G}_n$, since by the LIL result (23) for $\hat{G}_n$ and (24), we have $\|\hat{G}_n - G\|_\infty = O_{\mathrm{a.s.}}(\sqrt{n^{-1}\log\log n})$.

PROOF OF LEMMA 5.3. Considering the random variables $Z_i(\theta) = (G_\theta(X_i) - G(X_i))^2$ and using Lemma 5.1, we show that

$$
|K(\theta;\hat{G}_n) - K(\theta)| \leq c\|\hat{G}_n - G\|_\infty + \sup_{\theta\in\Theta}\left|\frac{1}{n}\sum_{i=1}^{n}(Z_i(\theta) - E(Z_i(\theta)))\right|,
$$

where $c$ is a nonnegative constant. The two terms on the right-hand sides no longer depend on $\theta$. The first tends to 0 with the desired rate of convergence by the LIL result given in (23). The second term is the supremum of an empirical process indexed by the functional class $\mathcal{H} = \{h(\cdot,\theta) = (G_\theta(\cdot) - G(\cdot))^2, \theta \in \Theta\}$ of Lipschitz bounded functions. Indeed, we have, by Lemma 5.2,

$$
\begin{aligned}
|h(x,\theta) - h(x,\theta')| &\leq |G_\theta(x) + G_{\theta'}(x) - 2G_\theta(x)| \times |G_\theta(x) - G_{\theta'}(x)| \\
&\leq c|\theta-\theta'|_2.
\end{aligned}
$$

Let $(\varepsilon_n)_{n\geq 1}$ be a sequence of real numbers decreasing to 0. It follows by a Bernstein type theorem of van der Vaart and Wellner ([41], page 246) that there exist nonnegative constants $A$ and $B$ such that

$$
P\left(\sup_{\theta\in\Theta}\left|\frac{1}{n}\sum_{i=1}^{n}(Z_i(\theta) - EZ_i(\theta))\right| > \varepsilon_n\right) \leq A(\sqrt{n}\varepsilon_n)^B\exp(-2n\varepsilon_n^2).
$$

It follows that if $\varepsilon_n = n^{-1/2+\alpha}$ with $\alpha > 0$, we get, by the Borel–Cantelli lemma,

$$
\sup_{\theta\in\Theta}\left|\frac{1}{n}\sum_{i=1}^{n}(Z_i(\theta) - EZ_i(\theta))\right| = o_{\mathrm{a.s.}}(n^{-1/2+\alpha}),
$$

which concludes the proof. $\square$

LEMMA 5.4. *Under* C1, *for* $k = 1, 2$, *there exists a real constant* $c > 0$ *such that, for all* $(\lambda_1, \lambda_2) \in [0, 1/2 - d]^2$ *and* $L \in L^\infty(\mathbb{R})$,

$$\left\| \left[ \frac{\partial^k}{\partial \lambda^k} A_\lambda S_r A_\lambda^{-1} \right]_{\lambda = \lambda_1} L - \left[ \frac{\partial^k}{\partial \lambda^k} A_\lambda S_r A_\lambda^{-1} \right]_{\lambda = \lambda_2} L \right\|_\infty \leq c |\lambda_1 - \lambda_2| \times \|L\|_\infty.$$

A straightforward consequence of the above lemma is that, for $k = 1, 2$ and $L \in L^\infty(\mathbb{R})$, there exists another real constant $c > 0$ such that

$$(31) \qquad \left\| \frac{\partial^k}{\partial \lambda^k} A_\lambda S_r A_\lambda^{-1} L \right\|_\infty \leq c \|L\|_\infty.$$

PROOF OF LEMMA 5.4. We prove the uniform Lipschitz property only for the case where $k = 1$, since the case where $k = 2$ uses the same technical arguments. For all $(\lambda_1, \lambda_2) \in [0, 1/2 - d]^2$, $L \in L^\infty(\mathbb{R})$, and all $x \in \mathbb{R}$, we have

$$\left| \left[ \frac{\partial}{\partial \lambda} A_\lambda S_r A_\lambda^{-1} \right]_{\lambda = \lambda_1} L(x) - \left[ \frac{\partial}{\partial \lambda} A_\lambda S_r A_\lambda^{-1} \right]_{\lambda = \lambda_2} L(x) \right|$$

$$\leq \left| \frac{1}{(1 - \lambda_1)^2} - \frac{1}{(1 - \lambda_2)^2} \right|$$

$$\times \sum_{k \geq 0} (k + 1) \left( \frac{\lambda_1}{1 - \lambda_1} \right)^k$$

$$(32) \qquad\qquad \times |L(-x + \mu_1 + \mu_2 + k\eta) - L(-x + 2\mu_2 + (k + 1)\eta)|$$

$$+ \frac{1}{(1 - \lambda_2)^2} \sum_{k \geq 0} (k + 1) \left| \left( \frac{-\lambda_1}{1 - \lambda_1} \right)^k - \left( \frac{-\lambda_2}{1 - \lambda_2} \right)^k \right|$$

$$\times |L(-x + \mu_1 + \mu_2 + k\eta) - L(-x + 2\mu_2 + (k + 1)\eta)|.$$

By the mean value theorem, there exist $\bar{\lambda}$ and $\tilde{\lambda}$ lying on the line segment with extremities $\lambda_1$ and $\lambda_2$ such that

$$\left| \frac{1}{(1 - \lambda_1)^2} - \frac{1}{(1 - \lambda_2)^2} \right| \leq \frac{2}{(1 - \bar{\lambda})^3} |\lambda_1 - \lambda_2|,$$

and for all $k \geq 0$,

$$\left| \left( \frac{-\lambda_1}{1 - \lambda_1} \right)^k - \left( \frac{-\lambda_2}{1 - \lambda_2} \right)^k \right| \leq k \left( \frac{\tilde{\lambda}}{1 - \tilde{\lambda}} \right)^{k-1} \left( \frac{1}{1 - \tilde{\lambda}} \right)^2 |\lambda_1 - \lambda_2|.$$

Using the above inequalities with (32), we obtain

$$\left\| \left[ \frac{\partial}{\partial \lambda} A_\lambda S_r A_\lambda^{-1} \right]_{\lambda = \lambda_1} L - \left[ \frac{\partial}{\partial \lambda} A_\lambda S_r A_\lambda^{-1} \right]_{\lambda = \lambda_2} L \right\|_\infty \leq \frac{12 \|L\|_\infty |\lambda_1 - \lambda_2|}{d^3},$$

which concludes the proof. $\square$

5.2. *Proof of Theorem* 2.1.

*Step* 1. Let $\{\sin(\alpha_1 t), \ldots, \sin(\alpha_p t)\}$ be a family of $p$ functions defined on $\mathbb{R}$. These functions are linearly independent if and only if we have

(33) $\quad \alpha_i \neq 0 \quad$ for $1 \leq i \leq p \quad$ and $\quad |\alpha_i| \neq |\alpha_j| \quad$ for $1 \leq i < j \leq p$.

Indeed, suppose that for $\beta_1, \ldots, \beta_p$ in $\mathbb{R}$ we have

$$\sum_{i=1}^{p} \beta_i \sin(\alpha_i t) = 0 \qquad \forall\, t \in \mathbb{R}.$$

Then, taking the derivative of the above expression with respect to $t$ at orders $1, 3, \ldots, 2p - 1$, we get at $t = 0$ the system of linear equations

$$\sum_{i=1}^{p} \beta_i \alpha_i^{2j+1} = 0 \qquad \text{for } 0 \leq j \leq p - 1.$$

The corresponding determinant is a Vandermonde type determinant different from 0 if and only if (33) is satisfied.

*Step* 2. We denote by $\Phi$ and $\Phi'$ the characteristic functions of $F$ and $F'$, respectively. Calculating the characteristic function of the two sides in (5), we get, for all $t \in \mathbb{R}$,

(34) $\quad (\lambda \exp(it\mu_1) + (1 - \lambda) \exp(it\mu_2))\Phi(t)$
$\qquad = (\lambda' \exp(it\mu_1') + (1 - \lambda') \exp(it\mu_2'))\Phi'(t).$

Since $F$ and $F'$ are c.d.f.'s of symmetric distributions, their characteristic functions are real continuous functions equal to 1 at $t = 0$. We have from (34) that the imaginary part of

$$(\lambda \exp(it\mu_1) + (1 - \lambda)\exp(it\mu_2))(\lambda' \exp(-it\mu_1') + (1 - \lambda')\exp(-it\mu_2'))$$

is equal to 0 in a neighborhood of 0. Then we have

(35) $\quad \lambda\lambda' \sin((\mu_1 - \mu_1')t) + \lambda(1 - \lambda')\sin((\mu_1 - \mu_2')t)$
$\qquad + (1 - \lambda)\lambda' \sin((\mu_2 - \mu_1')t) + (1 - \lambda)(1 - \lambda')\sin((\mu_2 - \mu_2')t) = 0$

on the whole real line, by analyticity of sine functions. We shall now consider two cases.

*Case* 1: $\lambda = 0$. Then (35) reduces to

(36) $\quad \lambda' \sin((\mu_2 - \mu_1')t) + (1 - \lambda')\sin((\mu_2 - \mu_2')t) = 0.$

If $\lambda' > 0$, then we have $1 - \lambda' > \lambda' > 0$, and by step 1 we need to consider the following cases:

- $\mu_2 = \mu_2'$ or $\mu_2 = \mu_1'$, hence by (36) $\mu_2' = \mu_1'$, which is not admissible.
- $|\mu_2 = \mu_1'| = |\mu_2 = \mu_2'|$, which by (36) leads to $\lambda' + (1 - \lambda') = 0$ (impossible) or $\lambda' - (1 - \lambda') = 0$ (not admissible).

It follows that $\lambda' = \lambda = 0$ and, hence, by (35) $\mu_2' = \mu_2$.

*Case* 2: $\lambda > 0$. From Case 1, we also have $\lambda' > 0$. Therefore, it remains to show that if $\mu_1 \neq \mu_2$, $\mu_1' \neq \mu_2'$, $(\lambda, \lambda') \in (0, 1/2)^2$ and that, for all $t \in \mathbb{R}$,

$$
\begin{aligned}
(37) \quad & \lambda\lambda' \sin((\mu_1 - \mu_1')t) + \lambda(1 - \lambda')\sin((\mu_1 - \mu_2')t) \\
& + (1 - \lambda)\lambda' \sin((\mu_2 - \mu_1')t) + (1 - \lambda)(1 - \lambda')\sin((\mu_2 - \mu_2')t) = 0,
\end{aligned}
$$

we have $(\lambda, \mu_1, \mu_2) = (\lambda', \mu_1', \mu_2')$. If we denote $\beta_1 = \lambda\lambda'$, $\beta_2 = \lambda(1 - \lambda')$, $\beta_3 = \lambda'(1 - \lambda)$ and $\beta_4 = (1 - \lambda)(1 - \lambda')$, then (37) is equivalent to

$$
\begin{aligned}
(38) \quad & \beta_1 \sin(\alpha t) + \beta_2 \sin((\alpha' - \eta)t) \\
& + \beta_3 \sin((\alpha + \eta)t) + \beta_4 \sin(\alpha' t) = 0 \qquad \forall\, t \in \mathbb{R},
\end{aligned}
$$

where $\alpha = \mu_1 - \mu_1'$, $\alpha' = \mu_2 - \mu_2'$ and $\eta = \mu_2 - \mu_1$. It is straightforward to see that if $\alpha = \alpha' = 0$, then $\lambda = \lambda'$. Then it remains to show that $(\alpha, \alpha') = (0, 0)$ is not admissible. To avoid a lengthy proof, we consider only the case $\alpha = 0$ and $\alpha' \neq 0$. The case $\alpha \neq 0$ and $\alpha' = 0$ is its symmetric counterpart and the case $\alpha \neq 0$ and $\alpha' \neq 0$ involves substantial calculations but is straightforward. Hence, if we suppose that $\alpha = 0$ and $\alpha' \neq 0$, equation (38) reduces to

$$
(39) \quad \beta_2 \sin((\alpha' - \eta)t) + \beta_3 \sin(\eta t) + \beta_4 \sin(\alpha' t) = 0 \qquad \forall\, t \in \mathbb{R}.
$$

Since $\alpha'$ and $\eta$ are nonnull, by Step 1, we have to consider the following cases:

- $\alpha' = \eta$: hence, $(\beta_3 + \beta_4)\sin(\eta t) = 0$ for all $t \in \mathbb{R}$. Then $\beta_3 + \beta_4 = 0$, which is not possible.
- $|\alpha' - \eta| = |\eta|$: hence, $\alpha' = 2\eta$. Then (39) reduces to

$$
(\beta_2 + \beta_3)\sin(\eta t) + \beta_4 \sin(2\eta t) = 0 \qquad \forall\, t \in \mathbb{R},
$$

  which, again by Step 1, cannot be satisfied for all $t \in \mathbb{R}$.
- Cases $|\alpha' - \eta| = |\alpha'|$ and $|\eta| = |\alpha'|$ lead respectively to $\alpha' = \eta/2$ and $\eta = -\alpha'$, hence, as in the previous case, the resulting equations cannot be satisfied for all $t \in \mathbb{R}$.

*Step* 3. Now, since $\lambda \in [0, 1/2)$ we have $|\lambda \exp(it\mu_1) + (1 - \lambda)\exp(it\mu_2)| \geq 1 - 2\lambda$. Then $\Phi = \Phi'$ and, finally, $F$ and $F'$ are equal.

### 5.3. *Proofs of Theorem* 3.1 *and Corollary* 3.1.

PROOF OF THEOREM 3.1. Let us write $\Phi_H$ for the characteristic function defined by $\Phi_H(t) = \int_{\mathbb{R}} \exp(itx)\,dH(x)$ for all $t \in \mathbb{R}$. Using the definitions of $A_\theta$ and $A_\theta^{-1}$ in (10), we obtain

$$
(40) \quad \Phi_{G_\theta}(t) = \frac{\lambda \exp(it\mu_1) + (1 - \lambda)\exp(it\mu_2)}{\lambda \exp(-it\mu_1) + (1 - \lambda)\exp(-it\mu_2)} \Phi_G(-t) \qquad \forall\, t \in \mathbb{R}.
$$

Moreover, because $\Phi_G(t) = (\lambda_0 \exp(it\mu_1^0) + (1 - \lambda_0) \exp(it\mu_2^0))\Phi_{F_0}(t)$ and $\Phi_{F_0}$ is an even function, $\Phi_{G_\theta}(t) = \Phi_G(t)$ for all $t \in \mathbb{R}$ implies that the imaginary part of

$$(\lambda \exp(it\mu_1) + (1 - \lambda) \exp(it\mu_2))(\lambda_0 \exp(-it\mu_1^0) + (1 - \lambda_0) \exp(-it\mu_2^0))$$

is null in a neighborhood of 0. Finally, by Step 2 of the proof of Theorem 2.1, we conclude that $\theta = \theta_0$.  □

PROOF OF COROLLARY 3.1.  Given the assumptions concerning $F_0$, we show that $G_\theta$ is a continuous function. By Theorem 3.1, if $\theta \neq \theta_0$, there exists $x_0 \in \mathbb{R}$ such that $G(x_0) \neq G_\theta(x_0)$, and there exist $\varepsilon > 0$ and $\alpha > 0$ such that $|G(x) - G_\theta(x)| > \varepsilon$ on $[x_0 - \alpha, x_0 + \alpha]$. It follows that

$$K(\theta) \geq \varepsilon^2 \int_{x_0 - \alpha}^{x_0 + \alpha} dG(x) = \varepsilon^2 (G(x_0 + \alpha) - G(x_0 - \alpha)) > 0.$$

Otherwise, if $\theta = \theta_0$ it is straightforward to check that $K(\theta) = 0$.  □

5.4.  *Proof of Theorem* 3.2.  Since the consistency proof for $\hat{\lambda}_n$ follows the lines of the consistency proof for $\hat{\theta}_n$ of problem (P2), it is omitted. For the remainder of this proof, we therefore suppose that $\hat{\lambda}_n$ converges almost surely to $\lambda_0$. By a first-order Taylor expansion of $\dot{K}(\cdot; \hat{G}_n)$ around $\hat{\lambda}_n$, we have

$$(41) \qquad \ddot{K}(\lambda_n^*; \hat{G}_n)\sqrt{n}(\hat{\lambda}_n - \lambda_0) = -\sqrt{n}\dot{K}(\lambda_0; \hat{G}_n),$$

where $\lambda_n^*$ lies on the line segment with extremities $\lambda_0$ and $\hat{\lambda}_n$. The desired result follows by proving the two statements

$$(42) \qquad\qquad \sqrt{n}\dot{K}(\lambda_0; \hat{G}_n) = O_P(1)$$

and

$$(43) \qquad\qquad \ddot{K}(\lambda_n^*; \hat{G}_n) \overset{\text{a.s.}}{\longrightarrow} 2 \int_{\mathbb{R}} \dot{G}_{\lambda_0}^2 \, dG > 0.$$

Result (42) follows from

$$|\dot{K}(\lambda_0; \hat{G}_n)| = \left| \int_{\mathbb{R}} 2\dot{\hat{G}}_{\lambda_0}^{(n)}(x)(\hat{G}_{\lambda_0}^{(n)}(x) - \hat{G}_n(x)) \, d\hat{G}_n(x) \right|$$

$$\leq 2\|\hat{G}_{\lambda_0}^{(n)} - \hat{G}_n\|_\infty \times \|\dot{\hat{G}}_{\lambda_0}^{(n)}\|_\infty$$

$$\leq 2\|A_{\lambda_0} S_r A_{\lambda_0}^{-1}[\hat{G}_n - G]\|_\infty \times \left\| \left[ \frac{\partial}{\partial \lambda} A_\lambda S_r A_\lambda^{-1} \right]_{\lambda = \lambda_0} \hat{G}_n \right\|_\infty.$$

The above inequality with Lemma 5.1 and (31) give the existence of a non-negative constant $c$ such that $|\dot{K}(\lambda_0; \hat{G}_n)| \leq c\|\hat{G}_n - G\|_\infty$. Thus, from result

(22), we get (42). In order to prove (43), let us write the second derivative of $K(\cdot; \hat{G}_n)$ at point $\lambda$,

$$\ddot{K}(\lambda; \hat{G}_n) = 2 \int_{\mathbb{R}} \ddot{G}_\lambda^{(n)} (\hat{G}_\lambda^{(n)} - \hat{G}_n) \, d\hat{G}_n + 2 \int_{\mathbb{R}} (\dot{G}_\lambda^{(n)})^2 \, d\hat{G}_n.$$

We have

(44)
$$\begin{aligned}
\left| \ddot{K}(\lambda_n^*; \hat{G}_n) - 2 \int_{\mathbb{R}} \dot{G}_{\lambda_0}^2 \, dG \right| \\
\leq |\ddot{K}(\lambda_n^*; \hat{G}_n) - \ddot{K}(\lambda_0; \hat{G}_n)| + \left| \ddot{K}(\lambda_0; \hat{G}_n) - 2 \int_{\mathbb{R}} \dot{G}_{\lambda_0}^2 \, dG \right|.
\end{aligned}$$

By very simple calculations, we show that

$$|\ddot{K}(\lambda_n^*; \hat{G}_n) - \ddot{K}(\lambda_0; \hat{G}_n)| \leq c|\lambda_n^* - \lambda_0|,$$

where $c$ is a nonnegative constant arising from Lemma 5.4, (31) and the fact that $\hat{G}_n$ is a cumulative distribution function. By the above inequality and the strong consistency of $\hat{\lambda}_n$, we conclude that $\ddot{K}(\lambda_n^*; \hat{G}_n) - \ddot{K}(\lambda_0; \hat{G}_n)$ converges almost surely to 0.

Concerning the second term of the right-hand side of (44), let us write

$$\begin{aligned}
\left| \ddot{K}(\lambda_0; \hat{G}_n) - 2 \int_{\mathbb{R}} \dot{G}_{\lambda_0}^2 \, dG \right| \\
\leq 2 \left| \int_{\mathbb{R}} \dot{G}_{\lambda_0}^2 \, d\hat{G}_n - \int_{\mathbb{R}} \dot{G}_{\lambda_0}^2 \, dG \right| + 2 \| \ddot{G}_{\lambda_0}^{(n)} \|_\infty \times \| \hat{G}_{\lambda_0}^{(n)} - \hat{G}_n \|_\infty \\
+ 2 (\| \dot{G}_{\lambda_0}^{(n)} \|_\infty + \| \dot{G}_{\lambda_0} \|_\infty) \times \| \dot{G}_{\lambda_0}^{(n)} - \dot{G}_{\lambda_0} \|_\infty.
\end{aligned}$$

Let us investigate the three terms on the right-hand side of the above inequality. From Lemmas 5.1 and 5.4, the second term is bounded, up to a multiplicative nonnegative constant, by

$$(\| \hat{G}_{\lambda_0}^{(n)} - G_{\lambda_0} \|_\infty + \| \hat{G}_n - G_{\lambda_0} \|_\infty) \leq (1 + \| A_{\lambda_0} S_r A_{\lambda_0}^{-1} \|) \times \| \hat{G}_n - G \|_\infty,$$

and then, tends to 0 almost surely, by using (21) and the LIL result given in (23). By similar arguments, we show that the third term has the same property. The first term is a centered empirical mean of i.i.d. random variables which, by Lemma 5.4, have a finite mean. Therefore, this term converges almost surely to 0 by the strong law of large numbers. Thus, it follows that

$$\left| \ddot{K}(\lambda_0; \hat{G}_n) - 2 \int_{\mathbb{R}} \dot{G}_{\lambda_0}^2 \, dG \right| = o_{\text{a.s.}}(1).$$

We conclude the proof, noticing that $\ddot{K}(\lambda_0) > 0$ [the proof, under C1, is similar to the proof of positive definitiveness of $\ddot{K}(\theta_0)$ in Section 5.5; therefore, it is omitted], and then

$$\ddot{K}(\lambda_0) = 2 \int_{\mathbb{R}} \dot{G}_{\lambda_0}^2 \, dG > 0.$$

5.5. *Proof of Theorem* 3.3.

*Proof of the consistency.* Our method is based on a consistency proof for miminum contrast estimators by Dacunha–Castelle and Duflo ([9], pages 94–96). Let us consider a countable dense set $D$ in $\Theta$. Then $\inf_{\theta \in \Theta} K_r(\theta; \hat{G}_n) = \inf_{\theta \in D} K_r(\theta; \hat{G}_n)$ is a measurable random variable. We define, in addition, the random variable

$$W(n, \xi) = \sup\{|K_r(\theta; \hat{G}_n) - K_r(\theta'; \hat{G}_n)|; (\theta, \theta') \in D^2, |\theta - \theta'|_2 \leq \xi\},$$

and recall that $K(\theta_0) = 0$. Let us consider a nonempty open ball $B_0$ centered on $\theta_0$ such that $K$ is bounded from below by a positive real number $2\varepsilon$ on $\Theta \backslash B_0$. Let us consider a sequence $(\xi_p)_{p \geq 1}$ decreasing to zero, and take $p$ such that there exists a covering of $\Theta \backslash B_0$ by a finite number $\ell$ of balls $(B_i)_{1 \leq i \leq \ell}$ with centers $\theta_i \in \Theta$, $i = 1, \ldots, \ell$, and radius less than $\xi_p$. Following Dacunha–Castelle and Duflo [9], we have

$$\begin{aligned}
(45) \quad \limsup_n \{\hat{\theta}_n \notin B_0\} &\subseteq \limsup_n \{W(n, \xi_p) > \varepsilon\} \\
&\cup \limsup_n \left\{\inf_{1 \leq i \leq \ell}(K_r(\theta_i; \hat{G}_n) - K_r(\theta_0; \hat{G}_n)) \leq \varepsilon\right\}.
\end{aligned}$$

By the uniform convergence result of Lemma 5.3, we have

$$(46) \quad P\left(\limsup_n \left\{\inf_{1 \leq i \leq \ell}(K_r(\theta_i; \hat{G}_n) - K_r(\theta_0; \hat{G}_n)) \leq \varepsilon\right\}\right) = 0.$$

Because $K$ is Lipschitz on $\Theta$ by Lemma 5.2, we have that, for sufficiently large $p'$, $|K(\theta) - K(\theta')| \leq \varepsilon/2$ for all $(\theta, \theta')$ such that $|\theta - \theta'|_2 \leq \xi_{p'}$. This implies

$$\begin{aligned}
\limsup_n &\{W(n, \xi_{p'}) > \varepsilon\} \\
&\subseteq \limsup_n \left\{2\sup_{\theta \in \Theta}|K_r(\theta; \hat{G}_n) - K(\theta)| + |K(\theta) - K(\theta')| > \varepsilon\right\} \\
&\subseteq \limsup_n \left\{2\sup_{\theta \in \Theta}|K_r(\theta; \hat{G}_n) - K(\theta)| > \varepsilon/2\right\},
\end{aligned}$$

and by Lemma 5.3 we have

$$P\left(\limsup_n \left\{2\sup_{\theta \in \Theta}|K_r(\theta; \hat{G}_n) - K(\theta)| > \varepsilon/2\right\}\right) = 0,$$

which leads to

$$(47) \quad P\left(\limsup_n \{W(n, \xi_{p'}) > \varepsilon\}\right) = 0.$$

By (45)–(47), we have proved the strong consistency of the contrast estimator $\hat{\theta}_n$.

*Proof for the convergence rate.* By standard Lebesgue theory, it is straightforward to show that, under C2, the contrast function $K$ is twice continuously differentiable on $\Theta$. If $\ddot{K}(\theta_0)$ is positive definite, by Corollary 3.1 and a Taylor expansion of $K$ of order 2 at $\theta_0$, there exist $\eta > 0$ and $\alpha > 0$ such that, for all $u$ satisfying $|u|_2 < \eta$ and $\theta_0 + u \in \overset{\circ}{\Theta}$,

$$(48) \qquad\qquad K(\theta_0 + u) \geq \alpha |u|_2^2.$$

For a column vector $v = (v_1, v_2, v_3)^T \in \mathbb{R}^3$, we have

$$(49) \qquad\qquad v^T \ddot{K}(\theta_0) v = 2 \int_{\mathbb{R}} (v^T \dot{G}_{\theta_0}(x))^2 \, dG(x) \geq 0.$$

By C2, we obtain that $x \mapsto \dot{G}_{\theta_0}(x)$ is continuous and that $G$ is continuous and strictly increasing on $\mathbb{R}$. Thus, (49) implies that $x \mapsto v^T \dot{G}_{\theta_0}(x)$ is the null function if $v^T \ddot{K}(\theta_0) v = 0$. Because under C2 we have $f_0' \in L^1(\mathbb{R})$, it is easy to show that $v^T \dot{g}_{\theta_0}(\cdot) \in L^1(\mathbb{R})$, where $g_\theta = A_\theta S_f A_\theta^{-1} g$. Moreover, using the Lebesgue derivation theorem and (40), and denoting $\eta_0 = \mu_2^0 - \mu_1^0$, we obtain

$$\begin{aligned}
\Phi_{v^T \dot{G}_{\theta_0}}(t) &= v^T \dot{\Phi}_{G_{\theta_0}}(t) \\
&= \frac{2\Phi_G(-t)}{(\lambda_0 \exp(-it\mu_1^0) + (1 - \lambda_0) \exp(-it\mu_2^0))^2} \\
&\quad \times [\cos(\eta_0 t)(v_1(1 - 2\lambda_0) + it(v_2(1 - \lambda_0) + v_3\lambda_0)) \\
&\qquad\qquad + it(v_2\lambda_0 + v_3(1 - \lambda_0))].
\end{aligned}$$

Therefore, $v^T \ddot{K}(\theta_0) v = 0$ implies that $\Phi_{v^T \dot{G}_{\theta_0}}$ is the null function. Because $\Phi_G(-t)/(\lambda_0 \exp(-it\mu_1^0) + (1 - \lambda_0) \exp(-it\mu_2^0))^2$ is not null in a neighborhood of 0, we obtain that the right multiplicative term of the right-hand side of the above equality is null in a neighborhood of 0, which in turn implies that $v = 0$. Thus $\ddot{K}(\theta_0)$ is positive definite and (48) holds.

Now, let us consider $B_0(\eta_n)$, the open ball centered at $\theta_0$ with radius $\eta_n > 0$. Notice that, for all $\theta \in \Theta \setminus B_0(\eta_n)$, we have $|\theta - \theta_0|_2 \geq \eta_n$. Then we write the event inclusions

$$\begin{aligned}
\{\hat{\theta}_n \notin B_0(\eta_n)\} &\subseteq \left\{ \inf_{\theta \in \Theta \setminus B_0(\eta_n)} K_r(\theta; \hat{G}_n) < K_r(\theta_0; \hat{G}_n) \right\} \\
&\subseteq \left\{ \inf_{\theta \in \Theta \setminus B_0(\eta_n)} K(\theta) - \sup_{\theta \in \Theta} |K_r(\theta; \hat{G}_n) - K(\theta)| < K_r(\theta_0; \hat{G}_n) \right\} \\
&\subseteq \left\{ \inf_{\theta \in \Theta \setminus B_0(\eta_n)} K(\theta) < 2 \sup_{\theta \in \Theta} |K_r(\theta; \hat{G}_n) - K(\theta)| \right\} \\
&\subseteq \left\{ \inf_{\theta \in \Theta \setminus B_0(\eta_n)} K(\theta) < \gamma_n \right\} \cup \left\{ \gamma_n \leq 2 \sup_{\theta \in \Theta} |K_r(\theta; \hat{G}_n) - K(\theta)| \right\}
\end{aligned}$$

for any arbitrary sequence $\gamma_n$. Thus, we have

$$\limsup_n \{\hat{\theta}_n \notin B_0(\eta_n)\} \subseteq \limsup_n \left\{ \inf_{\theta \in \Theta \setminus B_0(\eta_n)} K(\theta) < \gamma_n \right\}$$

$$\cup \limsup_n \left\{ \gamma_n \leq 2 \sup_{\theta \in \Theta} |K_r(\theta; \hat{G}_n) - K(\theta)| \right\}.$$

Choosing now $\gamma_n = n^{-1/2+\alpha}$ and $\eta_n = n^{-1/4+\beta/2}$, with $0 < \alpha < \beta$ taken arbitrarily small, it follows from (48) and the uniform almost sure rate of convergence of $K_r(\hat{G}_n)$ toward $K$, given in Lemma 5.3, that

$$P\left( \limsup_n \left\{ \inf_{\theta \in \Theta \setminus B_0(\eta_n)} K(\theta) < \gamma_n \right\} \right) = 0$$

and

$$P\left( \limsup_n \left\{ \gamma_n \leq 2 \sup_{\theta \in \Theta} |K_r(\theta; \hat{G}_n) - K(\theta)| \right\} \right) = 0.$$

In conclusion, $\hat{\theta}_n$ converges almost surely toward $\theta_0$ at rate $n^{-1/4+\delta}$, with $\delta > 0$ chosen arbitrarily small.

5.6. *Proof of Theorem* 3.4.

*Proof of* (i) *and* (ii). We have

$$\hat{F}_n - F_0 = \tfrac{1}{2}(I + S_r)[A_{\hat{\theta}_n}^{-1}\hat{G}_n - A_{\theta_0}^{-1}G].$$

Thus, there exists a nonnegative real constant $c$ such that

$$\|\hat{F}_n - F_0\|_\infty \leq \|A_{\hat{\theta}_n}^{-1}(\hat{G}_n - G)\|_\infty + \|(A_{\hat{\theta}_n}^{-1} - A_{\theta_0}^{-1})G\|_\infty$$

$$\leq \|A_{\hat{\theta}_n}^{-1}\| \times \|\hat{G}_n - G\|_\infty + c|\hat{\theta}_n - \theta_0|_2$$

$$\leq \frac{1}{2d}\|\hat{G}_n - G\|_\infty + c|\hat{\theta}_n - \theta_0|_2,$$

where the second inequality follows from (29) in the proof of Lemma 5.2 and the last inequality follows from (21), using the fact that $G$ is Lipschitz. Finally, the above inequality together with (22) [resp. (23)] and Theorem 3.2 (resp. Theorem 3.3) yield result (i) [resp. result (ii)].

*Proof of* (iii). By the Devroye [12] $L^1$-consistency result, we have

$$\|\hat{g}_n - g\|_1 = \int_{\mathbb{R}} |\hat{g}_n(x) - g(x)| \, dx \xrightarrow{\text{a.s.}} 0 \qquad (50)$$

as $n \to +\infty$, providing that $b_n \to 0$ and $nb_n \to +\infty$. Then we can write

$$\|\hat{f}_n - f_0\|_1 = \|A_{\hat{\theta}_n}^{-1}\hat{g}_n - A_{\theta_0}^{-1}g\|_1$$

$$\leq \|A_{\hat{\theta}_n}^{-1}(\hat{g}_n - g)\|_1 + \|(A_{\hat{\theta}_n}^{-1} - A_{\theta_0}^{-1})g\|_1$$

$$\leq \frac{1}{2d}\|\hat{g}_n - g\|_1 + C|\hat{\theta}_n - \theta_0|_2,$$

where the last inequality comes from (29), because $f_0' \in L^1(\mathbb{R})$ and, thus, the same holds for $g$. We conclude with Theorems 3.2 and 3.3, and (50).

L. Bordes
S. Mottelet
Université de Technologie de Compiègne
B.P. 529
60205 Compiègne cedex
France
E-mail: laurent.bordes@utc.fr
        stephane.mottelet@utc.fr

P. Vandekerkhove
Université de Marne-la-Vallée
Citée Descartes-5 Boulevard Descartes
Champs-sur-Marne
77454 Marne-la-Vallée cedex 2
France
E-mail: pierre.vandek@univ-mlv.fr